\documentclass[12pt,leqno]{article}

\newtheorem{theorem}{Theorem}
\newtheorem{lemma}[theorem]{Lemma}
\newtheorem{proposition}[theorem]{Proposition}
\newtheorem{definition}[theorem]{Definition}
\newtheorem{remark}[theorem]{Remark}

\newcommand{\begintheorem}{\addtocounter{equation}{1}\begin{theorem}}
\newcommand{\beginlemma}{\addtocounter{equation}{1}\begin{lemma}}
\newcommand{\beginproposition}{\addtocounter{equation}{1}\begin{proposition}}
\newcommand{\begindefinition}{\addtocounter{equation}{1}\begin{definition}}
\newcommand{\beginremark}{\addtocounter{equation}{1}\begin{remark}}

\begin{document}

\title{Potpourri, 3}

\author{Stephen William Semmes	\\
	Rice University		\\
	Houston, Texas}

\date{}

\maketitle


\renewcommand{\thefootnote}{}   

\footnotetext{These notes are connected to the ``potpourri''
topics class in the Department of Mathematics, Rice University,
in the fall semester of 2004.}

\tableofcontents

\bigskip

	Let $n$ be a positive integer, and consider the vector space
of real-valued continuous functions on the unit sphere in ${\bf R}^n$.
If $1 \le p < \infty$, then we can define the $L^p$ norm on this space
in the usual manner, by taking the $p$th power of the absolute value
of a function on the sphere, integrating it, and taking the $(1/p)$th
power of the result.  For $p = \infty$ we can use the supremum norm,
which is the maximum of the absolute value of a continuous function on
the sphere.

	Inside the space of continuous functions on the unit sphere in
${\bf R}^n$ we have the $n$-dimensional linear subspace of linear
functions, whose value at a point is given by the inner product of
that point with some fixed vector in ${\bf R}^n$.  For each $p$, $1
\le p \le \infty$, the $p$-norm mentioned in the previous paragraph
applied to a linear function reduces to a constant multiple of the
vector used in the inner product defining the linear function.

	Thus there is an $n$-dimensional inner product space embedded
isometrically into the space of real-valued continuous functions on
the unit sphere in ${\bf R}^n$ equipped with the $p$-norm for any $p$,
$1 \le p \le \infty$, and similarly one can get $n$-dimensional
complex inner product spaces in the space of complex-valued continuous
functions on the unit sphere in ${\bf C}^n$ equipped with the $p$-norm
for any $p$, $1 \le p \le \infty$.  In the first three sections of
these notes we shall consider some related constructions in which
$p$-norms can be compared to norms associated to inner products.

\section{Gaussian random variables}
\label{Gaussian random variables}
\setcounter{equation}{0}

	Recall that
\begin{equation}
	\int_{-\infty}^\infty \exp (- \pi x^2) \, dx = 1.
\end{equation}
Indeed, the integral is a positive real number whose square is equal
to
\begin{equation}
	\int_{{\bf R}^2} \exp (- \pi (x^2 + y^2)) \, dx \, dy,
\end{equation}
and this can be rewritten in polar coordinates as
\begin{equation}
   \int_0^{2 \pi} \int_0^\infty \exp (- \pi \, r^2) \, r \, dr \, d\theta
\end{equation}
which reduces to
\begin{equation}
	2 \pi \int_0^\infty r \exp( - \pi \, r^2) \, dr.
\end{equation}
This last integral can be computed directly, using the fundamental
theorem of calculus.

	Let $n$ be a positive integer, and let
\begin{equation}
	\langle v, w \rangle = \sum_{j=1}^n v_j \, w_j
\end{equation}
be the standard inner product on ${\bf R}^n$, $v = (v_1, \ldots,
v_n)$, $w = (w_1, \ldots, w_n)$.  Let $h(x)$ be a linear function on
${\bf R}^n$, so that there is a $v \in {\bf R}^n$ such that $h(x) =
\langle x, v \rangle$ for all $x \in {\bf R}^n$.

	It follows from the $1$-dimensional computation that
\begin{equation}
	\int_{{\bf R}^n} \exp (- \pi \langle x, x \rangle) \, dx = 1,
\end{equation}
i.e., the $n$-dimensional integral reduces exactly to a product of
$1$-dimensional integrals of the same type.  If $p$ is a positive real
number, consider
\begin{equation}
	\bigg(\int_{{\bf R}^n} |h(x)|^p \, 
		\exp (-\pi \langle x, x \rangle) \, dx \bigg)^{1/p}.
\end{equation}
This is equal to the product of $|v| = \langle v, v \rangle^{1/2}$ and
\begin{equation}
	\bigg(\int_{\bf R} |x|^p \, \exp (- \pi \, x^2) \, dx \bigg)^{1/p}.
\end{equation}
Of course this is trivial if $v = 0$, and if $v \ne 0$, then $h(x)$ is
constant in the directions orthogonal to $v$.  We can integrate out
those directions to get a $1$-dimensional integral, and after we pull
out $|v|$ the result depends only on $p$.

	Thus the linear functions on ${\bf R}^n$ form an
$n$-dimensional linear subspace of real-valued functions more
generally, e.g., real-valued continuous functions of polynomial
growth.  The $p$-norm reduces to a constant multiple of the Euclidean
norm for these functions.  This does not quite work for $p = \infty$,
in the sense that the functions are unbounded.  Of course the
coordinate functions $x_1, \ldots, x_n$ on ${\bf R}^n$, with respect
to the probability distribution given by the Gaussian function $\exp
(-\pi \, \langle x, x \rangle)$, are independent random variables with
mean $0$ and the same Gaussian distribution individually.  The linear
functions on ${\bf R}^n$ are the mean $0$ Gaussian random variables
which are linear combinations of these $n$ functions.

\section{Rademacher functions}
\label{Rademacher functions}
\setcounter{equation}{0}

	For each positive integer $\ell$ let $\mathcal{B}_\ell$ denote the
set of sequences $(x_1, \ldots, x_\ell)$ such that each $x_j$ is either
$-1$ or $1$.  Thus $\mathcal{B}_\ell$ has $2^\ell$ elements.

	We can think of $\mathcal{B}_\ell$ as a commutative group with
respect to coordinatewise multiplication.  In other words,
$\mathcal{B}_\ell$ is a product of $\ell$ copies of the group with $2$
elements.

	For each $j = 1, \ldots, \ell$ define $r_j$ to be the function
on $\mathcal{B}_\ell$ given by $r_j(x) = x_j$.  One can think of $r_j$
as a homomorphism from $\mathcal{B}_\ell$ into the group $\{\pm 1\}$.
One can also think of the $r_j$'s as random variables, with respect to
the uniform probability distribution on $\mathcal{B}_\ell$, so that
the $r_j$'s are independent random variables which have mean $0$ and
are identically distributed as fair coin tosses.

	If $I$ is a subset of $\{1, \ldots, \ell\}$, let us write
$w_I$ for the function on $\mathcal{B}_\ell$ which is the product of
the $r_j$'s with $j \in I$, where this is interpreted as being the
function which is equal to $1$ everywhere on $\mathcal{B}_\ell$ when
$I = \emptyset$.  The $r_j$'s are called Rademacher functions, and the
$w_I$'s are called Walsh functions.

	For a pair of real-valued functions $f_1(x)$, $f_2(x)$ on
$\mathcal{B}_\ell$, one can define their inner product to be
\begin{equation}
	2^{-\ell} \sum_{x \in \mathcal{B}_\ell} f_1(x) \, f_2(x).
\end{equation}
The inner product of a Walsh function with itself is equal to $1$,
since $w_I(x)^2 = 1$ for all subsets $I$ of $\{1, \ldots, \ell\}$ and
all $x \in \mathcal{B}_\ell$.  If $I_1$, $I_2$ are distinct subsets of
$\{1, \ldots, \ell\}$, then one can show that the inner product of the
corresponding Walsh functions on $\mathcal{B}_\ell$ is equal to $0$.

	Thus the Walsh functions form an orthonormal basis for the
vector space of real-valued functions on $\mathcal{B}_\ell$ with respect
to the inner product that we have defined.  That is, they are
orthonormal by the remarks in the previous paragraph, and they form a
basis because there are $2^\ell$ Walsh functions and the vector space of
functions on $\mathcal{B}_\ell$ has dimension equal to $2^\ell$.

	The Rademacher functions $r_1, \ldots, r_\ell$ are special
cases of Walsh functions, corresponding to subsets $I$ of $\{1,
\ldots, \ell\}$ with exactly one element.  The linear combinations of
the Rademacher functions form a very interesting $\ell$-dimensional
subspace of the vector space of real-valued functions on
$\mathcal{B}_\ell$.

	For $0 < p < \infty$, consider the quantity
\begin{equation}
	\bigg(2^{-\ell} \sum_{x \in \mathcal{B}_\ell} |f(x)|^p \bigg)^{1/p},
\end{equation}
where $f(x)$ is a real-valued function on $\mathcal{B}_\ell$.  When $p =
2$ this is the norm associated to the inner product on functions on
$\mathcal{B}_\ell$.  This quantity is monotone increasing in $p$, as a
result of the fact that $t^r$ defines a convex function on the
nonnegative real numbers when $r \ge 1$.

	For each $p \ge 2$ there is a positive real number
$C(p)$ such that
\begin{equation}
	\bigg(2^{-\ell} \sum_{x \in \mathcal{B}_\ell} |f(x)|^p \bigg)^{1/p}
		\le C(p) \,
	\bigg(2^{-\ell} \sum_{x \in \mathcal{B}_\ell} |f(x)|^2 \bigg)^{1/2}
\end{equation}
when $f$ is a linear combination of Rademacher functions.  To see
this, it suffices to consider the case where $p$ is an even integer,
by monotonicity in $p$.

	It is instructive to start with the case where $p = 4$.  If
$f(x)$ is a linear combination of Rademacher functions, then one can
expand $f(x)^4$ as a linear combination of certain Walsh functions,
many of which have sum equal to $0$, and the nonzero sums can be
estimated in terms of the $p = 2$ case.  Similar computations apply
for larger $p$'s.

	One can also show that for each positive real number $q \le 2$
there is a positive real number $C(p)$ such that
\begin{equation}
	\bigg(2^{-\ell} \sum_{x \in \mathcal{B}_\ell} |f(x)|^2 \bigg)^{1/2}
		\le C(p) \,
	\bigg(2^{-\ell} \sum_{x \in \mathcal{B}_\ell} |f(x)|^q \bigg)^{1/q}.
\end{equation}
One can derive this from the earlier result using H\"older's
inequality, i.e., the quantity for $p = 2$ can be estimated in terms
of the product of fractional powers of the quantity for $p = 4$ and
the quantity for $q < 2$ by H\"older's inequality, and therefore the
quantity for $p = 2$ can be estimated in terms of the quantity for $q
< 2$ since the quantity for $p = 4$ can be estimated in terms of the
quantity for $p = 2$.

\section{Lacunary series}
\label{lacunary series}
\setcounter{equation}{0}

	Suppose that $a_0, \ldots, a_n$ are complex numbers.
Because
\begin{equation}
	\frac{1}{2 \pi} \int_{\bf T} z^j \, \overline{z}^l \, |dz|
\end{equation}
is equal to $0$ when $j \ne l$ and equal to $1$ when $j = l$, where
${\bf T}$ denotes the unit circle in the complex plane, we have that
\begin{equation}
 \frac{1}{2 \pi} \int_{\bf T} \biggl| \sum_{j=0}^n a_j \, z^j \biggr|^2 \, |dz|
	= \sum_{j=0}^n |a_j|^2.
\end{equation}

	Now suppose that $c_0, \ldots, c_m$ are complex numbers, and
consider the function $f(z) = \sum_{j=0}^m c_j \, z^{2^j}$.  To
estimate
\begin{equation}
	\frac{1}{2 \pi} \int_{\bf T} |f(z)|^4 \, |dz|,
\end{equation}
one can write $|f(z)|^4$ as $f(z)^2 \, \overline{f(z)}^2$ and expand
the sums.  Many of the terms integrate to $0$.

	Indeed, the integral of
\begin{equation}
 z^{2^{j_1}} \, z^{2^{j_2}} \, \overline{z}^{2^{l_1}} \, \overline{z}^{2^{l_2}}
\end{equation}
over the unit circle is equal to $0$ unless
\begin{equation}
	2^{j_1} + 2^{j_2} = 2^{l_1} + 2^{l_2}.
\end{equation}
One can check that this happens for nonnegative integers $j_1, j_2,
l_1, l_2$ if and only if $j_1 = l_1$ and $j_2 = l_2$ or $j_1 = l_2$
and $j_2 = l_1$.

	As a result one can show that the integral of $|f(z)|^4$ over
the unit circle is bounded by a constant times the square of the
integral of $|f(z)|^2$.  There are many variants and extensions of
this observation.

\section{Some matrix norms}
\label{some matrix norms}
\setcounter{equation}{0}

	Fix positive integers $m$, $n$.  Suppose that $A$ is a linear
mapping from ${\bf R}^m$ to ${\bf R}^n$ associated to the $m \times n$
matrix $(a_{j, l})$ of real numbers.  Thus for $x = (x_1, \ldots, x_m)
\in {\bf R}^m$, the $l$th component of $A(x)$ is equal to
$\sum_{j=1}^m a_{j, l} \, x_j$.

	If $y = A(x)$ for some $x \in {\bf R}^m$, then
\begin{equation}
	|y_1| + \cdots + |y_n|
		\le \bigg(\sum_{j=1}^m \sum_{l=1}^n |a_{j, l}| \bigg)
			\, \max (|x_1|, \ldots, |x_m|).
\end{equation}
This inequality is optimal if the $a_{j, l}$'s are all nonnegative
real numbers, as one can see by taking $x_j = 1$ for all $j$.

	Now suppose that $T$ is a linear mapping from ${\bf R}^n$ to
${\bf R}^m$ associated to the $n \times m$ matrix $(t_{p, q})$ of real
numbers in the same way as before.  If $x = T(y)$, then
\begin{eqnarray}
\lefteqn{\max(|x_1|, \ldots, |x_m|)}		\\
	& & \le \max \{ |t_{p, q}| : 1 \le p \le n, 1 \le q \le m\}
			\ (|y_1| + \cdots + |y_n|),	\nonumber
\end{eqnarray}
and it is easy to see that this inequality is sharp.

	Let $\lambda_1, \ldots, \lambda_m$ denote the linear mappings
from ${\bf R}^m$ into ${\bf R}$ which take a given vector in ${\bf
R}^m$ to its $m$ coordinates.  Also let $e_1, \ldots, e_n$ denote the
standard basis vectors in ${\bf R}^n$, so that $e_j$ has $j$th
coordinate equal to $1$ and all other coordinates equal to $0$.  We
can express $A$ as a sum of rank-$1$ operators in the standard way,
\begin{equation}
	A(x) = \sum_{j=1}^m \sum_{l=1}^n a_{j, l} \, \lambda_l(x) \, e_j.
\end{equation}

	The composition of $T \circ A$ defines a linear mapping from
${\bf R}^m$ to itself, and the trace of the linear mapping is equal to
\begin{equation}
	\sum_{j=1}^m \sum_{l=1}^n a_{j, l} \, t_{l, j}.
\end{equation}
The absolute value of the trace of $T \circ A$ is less than or equal
to the product of $\sum_{j=1}^m \sum_{l=1}^n |a_{j, l}|$ and $\max
\{|t_{p, q}| : 1 \le p \le n, 1 \le q \le m\}$, and this inequality is
sharp.

\section{Grothendieck's inequality}
\label{grothendieck's inequality}
\setcounter{equation}{0}

	Let $m$, $n$ be positive integers, and let $(a_{j, l})$ be an
$m \times n$ matrix of real numbers.  Let us assume that this matrix
is \emph{restricted} in the sense that
\begin{equation}
	\biggl|\sum_{j=1}^m \sum_{l=1}^n a_{j, l} v_j \, w_l \biggr| \le 1
\end{equation}
whenever $v = (v_1, \ldots, v_m)$, $w = (w_1, \ldots, w_n)$ satisfy
\begin{equation}
	|v_1|, \ldots, |v_m|, |w_1|, \ldots, |w_n| \le 1.
\end{equation}
Equivalently, $(a_{j, l})$ is restricted if
\begin{equation}
	|y_1| + \cdots + |y_n| \le \max (|x_1|, \ldots, |x_m|)
\end{equation}
for all $x \in {\bf R}^m$, where $y \in {\bf R}^n$ is given by $y_l =
\sum_{j=1}^m a_{j, l} \, x_j$.

	Let $V$ be a finite-dimensional real vector space equipped
with an inner product $\langle v, w \rangle$, which one can simply
take to be a Euclidean space with the standard inner product.  Let
$v_1, \ldots, v_m$ and $w_1, \ldots, w_n$ be vectors in $V$.
Grothendieck's inequality states that there is a universal constant $k
> 0$ so that if
\begin{equation}
	\|v_1\|, \ldots, \|v_m\|, \|w_1\|, \ldots, \|w_n\| \le 1,
\end{equation}
where $\|u\| = \langle u, u \rangle^{1/2}$ is the norm of $u \in V$
associated to the inner product, then
\begin{equation}
 \biggl|\sum_{j=1}^m \sum_{l=1}^n a_{j, l} \langle v_j, w_l \rangle \biggr|
	\le k.
\end{equation}
See volume 1 of \cite{L-T2} for a proof, with $k = (\exp (\pi/2) -
\exp (-\pi/2))/2$.

\section{Maximal functions}
\label{maximal functions}
\setcounter{equation}{0}

	As in Section \ref{Rademacher functions}, for each positive
integer $\ell$ let $\mathcal{B}_\ell$ denote the set of sequences $x =
(x_1, \ldots, x_\ell)$ of length $\ell$ such that $x_j \in \{\pm 1\}$
for all $j$.  In this section and the next one, for each real-valued
function $f$ on $\mathcal{B}_\ell$, we put
\begin{equation}
	\|f\|_p = 
   \bigg( 2^{-\ell} \sum_{x \in \mathcal{B}_\ell} |f(x)|^p \bigg)^{1/p}
\end{equation}
when $1 \le p < \infty$ and
\begin{equation}
	\|f\|_\infty = \max \{|f(x)| : x \in \mathcal{B}_\ell\}.
\end{equation}

	If $x \in \mathcal{B}_\ell$ and $k$ is an integer with $0 \le
k \le \ell$, let $N_k(x)$ denote the set of $y \in \mathcal{B}_\ell$
such that $y_j = x_j$ when $j \le k$.  Thus $N_0(x) =
\mathcal{B}_\ell$, $N_\ell(x) = \{x\}$, and $N_k(x)$ has $2^{\ell -
k}$ elements for each $k$.

	If $f(x)$ is a real-valued function on $\mathcal{B}_\ell$ and
$0 \le k \le \ell$, define $E_k(f)$ to be the function on
$\mathcal{B}_\ell$ whose value at a given point $x$ is the average of
$f$ on $N_k(x)$, i.e.,
\begin{equation}
	E_k(f)(x) = 2^{-\ell + k} \sum_{y \in N_k(x)} f(y).
\end{equation}
Thus $E_k(f)(x)$ is a function on $\mathcal{B}_\ell$ which really only
depends on the first $k$ coordinates of $x$.  If $h(x)$ is a function
on $\mathcal{B}_\ell$ which depends on only the first $k$ coordinates
of $x$, then
\begin{equation}
	E_k(h) = h,
\end{equation}
and in fact
\begin{equation}
	E_k(f \, h) = h \, E_k(f)
\end{equation}
for all functions $f$ on $\mathcal{B}_\ell$.

	Notice that if $0 \le k < \ell$ and $f$ is a function on
$\mathcal{B}_\ell$ such that $f(x)$ depends only on the first $k + 1$
coordinates of $x$, then $f$ can be expressed uniquely as $f_1 +
r_{k+1} \, f_2$, where $r_{k+1}(x) = x_{k+1}$ as in Section
\ref{Rademacher functions} and $f_1(x)$, $f_2(x)$ depend only on the
first $k$ coordinates of $x$.  Moreover, $E_k(f) = f_1$ in this
situation.

	For each real-valued function $f$ on $\mathcal{B}_\ell$
the associated maximal function is given by
\begin{equation}
	M(f)(x) = \max \{|E_k(x)| : 0 \le k \le \ell\},
\end{equation}
$x \in \mathcal{B}_\ell$.  Observe that
\begin{equation}
	\|M(f)\|_\infty \le \|f\|_\infty
\end{equation}
for all functions $f$ on $\mathcal{B}_\ell$, and in fact we have
equality, since $E_\ell(f) = f$.

	Fix a function $f$ on $\mathcal{B}_\ell$.  Let $\lambda$ be a
positive real number, and consider the set
\begin{equation}
	A_\lambda = \{x \in \mathcal{B}_\ell : M(f) > \lambda\}.
\end{equation}
Thus $A_\lambda$ is the set of $x \in \mathcal{B}_\ell$ such that
$|E_k(f)(x)| > \lambda$ for some $k$, $0 \le k \le \ell$.

	If $x \in A_\lambda$, so that $|E_k(f)(x)| > \lambda$ for some
$k$, $0 \le k \le \ell$, then we automatically have that $|E_k(f)(y)|
> \lambda$ for all $y \in N_k(x)$, since $E_k(f)$ is constant on
$N_k(x)$.  We can describe $A_\lambda$ as the union of the subsets of
$\mathcal{B}_\ell$ of the form $N_k(x)$ on which the absolute value of
the average of $f$ is larger than $\lambda$.

	If $x$, $x'$ are elements of $\mathcal{B}_\ell$ and $k$, $k'$
are integers such that $0 \le k, k' \le \ell$, then
\begin{equation}
	N_k(x) \subseteq N_{k'}(x'), \hbox{ or }
	N_{k'}(x') \subseteq N_k(x), \hbox{ or }
	N_k(x) \cap N_{k'}(x') = \emptyset.
\end{equation}
This is easy to check from the definitions.  We can think of
$A_\lambda$ as the union of the maximal subsets of $\mathcal{B}_\ell$
of the form $N_k(x)$ on which the absolute value of the average of $f$
is larger than $\lambda$, and the maximality of these subsets implies
that they are pairwise disjoint.

	Let us write $|A_\lambda|$ for $2^{-\ell}$ times the number of
elements of $A_\lambda$.  In other words, if $a_\lambda(x)$ is the
function on $\mathcal{B}_\ell$ such that $a_\lambda(x) = \lambda$ when
$x \in A_\lambda$ and $a_\lambda(x) = 0$ otherwise, then
\begin{equation}
	\|a_\lambda\|_1 = \lambda \, |A_\lambda|.
\end{equation}

	A key point now is that
\begin{equation}
	\lambda \, |A_\lambda| \le \|f\|_1.
\end{equation}
This follows from the fact that $A_\lambda$ is a disjoint union of
sets of the form $N_k(x)$ on which the absolute value of the average
of $f$ is larger than $\lambda$.

	Let $f'(x)$ be the function on $\mathcal{B}_\ell$ which is
equal to $f(x)$ when $|f(x)| > \lambda$ and which is equal to $0$
otherwise.  Of course
\begin{equation}
	|f(x) - f'(x)| \le \lambda
\end{equation}
for all $x \in \mathcal{B}_\ell$.  Hence
\begin{equation}
	M(f)(x) \le M(f')(x) + \lambda
\end{equation}
for all $x \in \mathcal{B}_\ell$.

	Therefore $A_{2 \lambda}$ is contained in the set of $x \in
\mathcal{B}_\ell$ such that $M(f') > \lambda$, and the previous
estimate implies that
\begin{equation}
	\lambda \, |A_{2 \lambda}| \le \|f'\|_1.
\end{equation}
Using this one can show that for each $p > 1$, $\|M(f)\|_p$ is bounded
by a constant depending only on $p$ times $\|f\|_p$.

\section{Square functions}
\label{square functions}
\setcounter{equation}{0}

	If $f_1$, $f_2$ are real-valued functions on
$\mathcal{B}_\ell$, let us write $(f_1, f_2)$ for their inner product,
\begin{equation}
	(f_1, f_2) = 2^{-\ell} \sum_{x \in \mathcal{B}_\ell} f_1(x) \, f_2(x).
\end{equation}
Thus the norm $\|f\|_2$ is the same as the norm associated to the
inner product, which is to say that $\|f\|_2^2 = (f, f)$ for all
real-valued functions $f$ on $\mathcal{B}_\ell$.

	For any function $f$ on $\mathcal{B}_\ell$, we can decompose
$f$ into the sum
\begin{equation}
	f = E_0(f) + \sum_{k=1}^\ell (E_k(f) - E_{k-1}(f)).
\end{equation}
The functions $E_0(f)$, $E_k(f) - E_{k-1}(f)$, $1 \le k \le \ell$, are
pairwise orthogonal with respect to the inner product just defined.

	In particular,
\begin{equation}
  \|f\|_2^2 = \|E_0(f)\|_2^2 + \sum_{k=1}^\ell \|E_k(f) - E_{k-1}(f)\|_2^2,
\end{equation}
where actually $E_0(f)$ is a constant and its norm reduces to its
absolute value.  Define the square function associated to $f$ by
\begin{equation}
	S(f)(x) = 
 \bigg(|E_0(f)(x)|^2 + 
		\sum_{k=1}^\ell |E_k(f)(x) - E_{k-1}(f)(x)|^2 \bigg)^{1/2}.
\end{equation}
The previous formula for $\|f\|_2$ can be rewritten as
\begin{equation}
	\|S(f)\|_2 = \|f\|_2.
\end{equation}

	As a special case, suppose that $f = \sum_{j=1}^\ell a_j \,
r_j$, where $a_1, \ldots, a_\ell$ are real numbers and $r_j(x) = x_j$.
Then $E_0(f) = 0$, $E_k(f) - E_{k-1}(f) = a_k \, r_k$, and $S(f)$ is
the constant $\big(\sum_{j=1}^\ell a_j^2 \big)^{1/2}$.

	Now let $f$ be any real-valued function on $\mathcal{B}_\ell$.
A very cool fact is that $\|S(f)\|_4$ is bounded by a constant times
$\|f\|_4$.

	It is easy enough to write $S(f)^4$ explicitly, by multiplying
out the sums.  One can rewrite this as a sum over $k$, where each part
in the sum is a product of something around level $k$ times a sum
involving the next levels.  When one sums over $x \in
\mathcal{B}_\ell$, these sums over levels $\ge k$ can be analyzed
using orthogonality.

	The basic conclusion is that $\|S(f)\|_4$ can be estimated,
using the Cauchy--Schwarz inequality, in terms of the product of
$\|S(f)\|_4^{1/2}$ and $\|M(f^2)\|_2$.  Because we can estimate
$\|M(f^2)\|_2$ in terms of $\|f^2\|_2$, which is the same as
$\|f\|_4^{1/2}$, we can estimate $\|S(f)\|_4$ in terms of $\|f\|_4$,
as desired.

\end{document}